\documentclass{amsart}
\usepackage{graphicx}

\def\Date{\Month\ \the\day, \the\year}  
\def\Month{\ifcase\number\month \relax\or January\or February\or
        March\or April\or May\or June\or July\or August\or September\or
        October\or November\or December\else\relax\fi}

\newtheorem{theorem}{Theorem}

\newtheorem{lemma}[theorem]{Lemma}
\theoremstyle{definition}
\newtheorem*{remark}{Remark}
\newtheorem*{example}{Theorem}

\long\def\ignore#1{}

\begin{document}

\title{On consecutive quadratic non-residues: a conjecture of Issai Schur}
\author{Patrick Hummel}
\address{1170 E. Del Mar Blvd., Apt. 10
        Pasadena, CA 91106-3425}
\email{hummel@its.caltech.edu}
\maketitle

Issai Schur once asked if it was possible to determine a bound,
preferably using elementary methods, such that for all prime
numbers $p$ greater than the bound, the greatest possible number
of consecutive quadratic non-residues modulo $p$ is always less
than $p^{1/2}$. (One can find a brief discussion of this problem
in R.~K.~Guy's book \cite{GUY}).  Schur also pointed out that the
greatest number of consecutive quadratic non-residues exceeds
$p^{1/2}$ for $p=13$, since $5$, $6$, $7$, and $8$ are all
quadratic non-residues (mod $p$).  This paper uses elementary
methods to prove the following:

\begin{example}  {\em $p=13$ is the only prime number for which the greatest
number of consecutive quadratic non-residues modulo $p$ exceeds
$p^{1/2}$.}

\end{example}

This problem has been attacked previously using both analytic and
elementary methods.  We shall briefly consider the results given
to us by analytic number theory, and then focus on the elementary
methods for the remainder of the paper.

In \cite{BU}, D.~A.~Burgess proves the following:

\begin{example}[D.~A.~Burgess \cite{BU}]  {\em If $\chi$
is any non-trivial Dirichlet character of prime modulus $p$ and
$\chi(N+1)=\chi(N+2)=...=\chi(N+H)$, then $H=O(p^{1/4}$} log {\em
$p)$.}

\end{example}

From this, it follows that there must be some $M$, such that for
all $p>M$, the greatest number of consecutive quadratic
non-residues modulo $p$ is less than $p^{1/2}$.  In \cite{KKN},
K.~K.~Norton asserts that he can refine Burgess's method to obtain
the following result:

\begin{example}[K.~K.~Norton \cite{KKN}]  {\em In Burgess's
theorem, $H<4.1p^{1/4}$} log {\em $p$ for all $p$.  If
$p>e^{15}\approx3.27\times10^{6}$, then $H<2.5p^{1/4}$} log {\em
$p$.}

\end{example}

This result implies that $M=e^{15}\approx3.27\times10^{6}$ is a
suitable value for the aforementioned constant. Unfortunately,
however, Norton does not prove this result in his paper, and
without a value such as $4.1$ for the implied constant in
Burgess's thoerem, we cannot use Burgess's theorem to find a
suitable constant, $M$, in Schur's conjecture.

\bigskip

Now we consider the work that has been done on the problem using
elementary methods.  A.~Brauer \cite{BIII} has proved the
following theorem:

\begin{example}[A.~Brauer \cite{BIII}]  {\em For prime numbers $p$, of the form $4n-1$,
the maximum length $l$, of sequences of quadratic residues and
non-residues satisfies $l<p^{1/2}$.}

\end{example}

R.~H.~Hudson then considers the case $p\equiv{1}$ (mod $4$) by
breaking it up into several cases.  In \cite{HII}, he proves that
the maximum number of consecutive quadratic non-residues (mod $p$)
is less than $p^{1/2}$ if $p\equiv{1}$ (mod $24$).  In \cite{H},
he demonstrates that this also holds if $p\equiv{5}$ or
$p\equiv{17}$ (mod $24$).  Putting these together, Hudson obtains
the following beneficial result:

\begin{example}[R.~H.~Hudson \cite{H}]  {\em If $p$ is a prime, and the greatest
number of consecutive quadratic non-residues modulo $p$ exceeds
$p^{1/2}$, then $p\equiv{13}$ (mod $24$).}

\end{example}

In the same paper, Hudson proposes a proof that the greatest
number of consecutive quadratic non-residues modulo $p$, is less
than $p^{1/2}$ for $p>2^{332}$ and $p\equiv{13}$ (mod $24$).  If
his proof of this result were correct, it would complete an
elementary proof of Schur's conjecture.  But upon a careful
reading, one sees that the proof of his assertion is flawed.  In
particular, Hudson claims that the existence of a quadratic
non-residue in the interval $(\frac{p^{1/2}}{128}-2^{3/2}p^{1/8},
\frac{p^{1/2}}{128})$ implies a quadratic {\em non-residue} is
contained in the interval $(p^{1/2}-2^{17/2}p^{1/8}, p^{1/2})$,
whereas it really implies a quadratic {\em residue} is contained
in the interval $(p^{1/2}-2^{17/2}p^{1/8}, p^{1/2})$, since $2$ is
a quadratic non-residue modulo $p$ if $p\equiv{13}$ (mod $24$).

\bigskip

The research presented here tackles the case $p\equiv{13}$ (mod
$24$) by using an alternate method.  In addition to completing an
elementary proof of Schur's conjecture, this paper also shows that
$p=13$ is the only prime number for which the greatest number of
consecutive quadratic non-residues (mod $p$) exceeds $p^{1/2}$.

The argument in this paper breaks down into two parts.  In the
first part, it is shown that if $p\equiv 13$ (mod $24$) is
sufficiently large and if there exists an interval containing more
than $p^{1/2}$ integers, all of which are quadratic non-residues
modulo $p$, then there must exist such an interval $J$ satisfying

\[
J\subset{(\frac{p+3+p^{1/2}}{2},
\frac{p}{2}+2^{1/2}p^{3/4}-p^{1/2})}.
\]

Then, in the second part, it is shown that this cannot hold for
$p$ sufficiently large by demonstrating that one can find two
integers, $a,b \in J$ such that if $R\equiv ab$ (mod $p$) and
$0\leq R\leq p-1$, then $R\in J$.  This implies that $J$ contains
a quadratic residue $R$, so $J$ could not have existed in the
first place.  This implies that the greatest number of consecutive
quadratic non-residues modulo $p$, is less than $p^{1/2}$ when $p$
is sufficiently large.  The remaining cases, where $p$ is less
than a given bound, are handled by computer.

This research is the result of work done under the mentorship of
Prof.~Dinakar Ramakrishnan while on a Caltech Summer Undergraduate
Research Fellowship.  I am deeply grateful to Prof.~Ramakrishnan
for his guidance and to David Whitehouse for reviewing my paper.

\begin{lemma}  If $p$ is a prime number such that
$p\equiv{13}$ (mod $24$), $p>38659$, and there is a sequence of
more than $p^{1/2}$ consecutive quadratic non-residues (mod $p$),
there must be such a sequence in the interval
$(\frac{p+3+p^{1/2}}{2}, \frac{p}{2}+2^{1/2}p^{3/4}-p^{1/2})$.

\end{lemma}

\proof

$p\equiv{13}$ (mod $24$) implies that every number of the form
$2a^{2}$ is a quadratic non-residue.  This means that there must
exist a quadratic non-residue, say $N$, in the interval

\[
(p^{1/2}-2^{3/2}p^{1/4}+2, p^{1/2}),
\]
because if $c$ is the smallest positive integer such that
$2c^{2}>p^{1/2}$, then $2(c-1)^{2}>p^{1/2}-2^{3/2}p^{1/4}+2$.

Suppose $J$ is an integer interval containing more than $p^{1/2}$
consecutive quadratic non-residues.  Multiplying each member of
$J$ by $N$ and reducing (mod $p$), we obtain a collection of
quadratic residues in which each quadratic residue differs from
the next by $N<p^{1/2}$.  This collection must span more than

\[
(p^{1/2}-1)(p^{1/2}-2^{3/2}p^{1/4}+2)>p-2^{3/2}p^{3/4}+p^{1/2}
\]
integers.  Since $-1$ is a quadratic residue (mod $p$), $p-b$ must
be a quadratic non-residue whenever $b$ is.  Therefore, if our
collection of quadratic residues is to lie entirely outside a
sequence of more than $p^{1/2}$ consecutive quadratic
non-residues, $J$ must either be contained in

\[
(1, 2^{1/2}p^{3/4}-\frac{p^{1/2}}{2}),
\]
\[
(p-2^{1/2}p^{3/4}+\frac{p^{1/2}}{2}, p-1),
\]
or
\begin{equation}\label{Equation 19}
(\frac{p}{2}-2^{1/2}p^{3/4}+\frac{p^{1/2}}{2},
\frac{p}{2}+2^{1/2}p^{3/4}-\frac{p^{1/2}}{2}).
\end{equation}
But $J$ cannot be fully contained in the first of these intervals
because the difference between the square numbers in $(1,
2^{1/2}p^{3/4}-\frac{p^{1/2}}{2})$ is less than $p^{1/2}$.
Similarly, $J$ cannot be contained in the second of these
intervals because any sequence of the form given by $J$ in
$(p-2^{1/2}p^{3/4}+\frac{p^{1/2}}{2}, p-1)$ would have to
correspond to a similar sequence in $(1,
2^{1/2}p^{3/4}-\frac{p^{1/2}}{2})$ So such a $J$ can only be
contained in the interval given by (\ref{Equation 19}).

We now refer to the following theorem of A. Brauer's:

\begin{example}[A. Brauer \cite{B}]  {\em The least odd quadratic
non-residue $u$ modulo a prime $p$ satisfies
$u<2^{3/5}p^{2/5}+2^{-(6/5)}\cdot25p^{1/5}+3$ for $p=8n+5$.}

\end{example}

This implies that there exists an odd quadratic non-residue $u$,
less than $p^{1/2}$ if $p>38659$.  Then, since $\frac{p+1}{2}$ is
a quadratic non-residue (mod $p$),
$u(\frac{p+1}{2})\equiv{\frac{p+u}{2}}$ is a quadratic residue
(mod $p$). Therefore, there exists a quadratic residue in the
interval ($\frac{p}{2}$, $\frac{p+p^{1/2}}{2}$), so there must
exist a corresponding quadratic residue in the interval
($\frac{p-p^{1/2}}{2}$, $\frac{p}{2}$), which means that $J$
cannot pass through $\frac{p}{2}$.

Combining this with the fact that $-b$ is a quadratic non-residue
whenever $b$ is, we know that if such a $J$ exists, there must be
at least one such $J$ in the interval

\begin{equation}\label{Equation 16}
(\frac{p}{2}, \frac{p}{2}+2^{1/2}p^{3/4}-\frac{p^{1/2}}{2}).
\end{equation}

Now note that $3$ is a quadratic residue (mod $p$).  Therefore,
for odd $m$, $\frac{p+3m}{2}$ must be a quadratic residue (mod
$p$) if $\frac{p+m}{2}$ is. Combined with the fact that there
exists a quadratic residue in the interval ($\frac{p}{2}$,
$\frac{p+p^{1/2}}{2}$), we find that if such a $J$ lies in the
interval given by (\ref{Equation 16}), that same $J$ must also lie
in the interval

\begin{equation}\label{Equation 17}
\left(\frac{p+3+p^{1/2}}{2},
\frac{p}{2}+2^{1/2}p^{3/4}-\frac{p^{1/2}}{2}\right).
\end{equation}

To see why, suppose that $\frac{p+3}{2}+x$ is the first entry in
$J$.  Then, we can assume that $\frac{p+1}{2}+x$ is a quadratic
residue, meaning $\frac{p}{2}+3(\frac{1}{2}+x)=\frac{p+3}{2}+3x$
is also a quadratic residue.  Therefore we must have
$\frac{p+3}{2}+3x-(\frac{p+3}{2}+x)>p^{1/2}$, or
$x>\frac{p^{1/2}}{2}$.

\endproof

\begin{lemma}\label{Lemma 2}  Suppose $p>38659$ is a prime congruent to $13$ modulo $24$ and $\frac{p+1}{2}+k$
is a quadratic non-residue, where $k>0$ is some fixed integer.
Then, if there exists $a$ such that
$\frac{1}{4}\leq{a}\leq\frac{15}{32}$, and
$(ap^{1/2}-2)^{2}>2k+2(1-a)p^{1/2}+2-\lfloor$$p^{1/2}$$\rfloor$,
and the difference between
$(\frac{p+1}{2})^{2}+k+k^{2}+2kx+x+x^{2}$ at
$x=\lfloor$$(1-a)p^{1/2}$$\rfloor$ and
$x=\lfloor$$ap^{1/2}$$\rfloor$$-2$ is greater than $p$,
$\frac{p+1}{2}+k$ is not the smallest number in a sequence of more
than $p^{1/2}$ consecutive quadratic non-residues.

\end{lemma}

\proof

Suppose that all the integers of the form

\[
\frac{p+1}{2}+k+m
\]
are quadratic non-residues, where $m$ is an integer ranging from
$0$ to $\lfloor$$p^{1/2}$$\rfloor$.

Note that the product of two integers of this form, say
$\frac{p+1}{2}+k+m$ and $\frac{p+1}{2}+k+n$, is a quadratic
residue, and equals

\begin{equation}\label{Equation 5}
\left(\frac{p+1}{2}\right)^{2}+(p+1)k+k^{2}+(m+n)(k+\frac{p+1}{2})+mn
\end{equation}
which if $m$ and $n$ both equal the same value, say $x$, reduces
to

\[
\left(\frac{p+1}{2}\right)^{2}+(p+1)k+k^{2}+(2x)(k+\frac{p+1}{2})+x^{2}
\]
which is congruent to

\begin{equation}\label{Equation 7}
\left(\frac{p+1}{2}\right)^{2}+k+k^{2}+2kx+x+x^{2}
\end{equation}
modulo $p$.

If there exists an $a$ such that
$\frac{1}{4}\leq{a}\leq\frac{15}{32}$ and the difference between
$(\frac{p+1}{2})^{2}+k+k^{2}+2kx+x+x^{2}$ at
$x=\lfloor$$(1-a)p^{1/2}$$\rfloor$ and
$x=\lfloor$$ap^{1/2}$$\rfloor$$-2$ is greater than $p$, we note
that we can find an integer $x$ contained in the interval
$(ap^{1/2}-2, (1-a)p^{1/2}]$, and an integer $c$ such that

\[
\left(\frac{p+1}{2}\right)^{2}+k+k^{2}+2kx+x+x^{2}>(c+\frac{1}{2})p+\frac{1}{2}+k+\lfloor
p^{1/2}\rfloor
\]
and

\begin{equation}\label{Equation 9}
\left(\frac{p+1}{2}\right)^{2}+k+k^{2}+2k(x-1)+x-1+(x-1)^{2}\leq(c+\frac{1}{2})p+\frac{1}{2}+k+\lfloor
p^{1/2}\rfloor.
\end{equation}

Now suppose that

\[
\left(\frac{p+1}{2}\right)^{2}+k+k^{2}+2k(x-1)+x-1+(x-1)^{2}\geq(c+\frac{1}{2})p+\frac{1}{2}+k.
\]
Then, combining this with (\ref{Equation 9}), we reach the absurd
conclusion that a quadratic residue equals a quadratic
non-residue. Therefore, we have

\begin{equation}\label{Equation 13}
\left(\frac{p+1}{2}\right)^{2}+k+k^{2}+2k(x-1)+x-1+(x-1)^{2}<(c+\frac{1}{2})p+\frac{1}{2}+k.
\end{equation}

Now consider (\ref{Equation 5}) again.  Let $m$ and $n$ vary so
that $m=x-y$ and $n=x+y$, where $x$ is an integer that satisfies
the above conditions, and $y$ is an integer ranging from $0$ to
the smallest integer larger than $ap^{1/2}-2$. Since $x$ lies in
$(ap^{1/2}-2$, $(1-a)p^{1/2}]$, we continue to meet the condition
that $m$ and $n$ are both integers between $0$ and
$\lfloor$$p^{1/2}$$\rfloor$ inclusive, because $y<ap^{1/2}-1$, and
$\lfloor$$p^{1/2}$$\rfloor$$-(1-a)p^{1/2}>ap^{1/2}-1$. If $m$ and
$n$ vary this way, the only part of (\ref{Equation 5}) that
changes is the product $mn$. Also note that

\[
0<mn-(m-1)(n+1)=n-m+1<p^{1/2}
\]
when $0\leq{n-m}<p^{1/2}-1$, which holds when $y$ varies as above.
So we have a collection of quadratic residues in which no
quadratic residue exceeds the next by more than $p^{1/2}$.  This
collection spans an interval of

\begin{equation}\label{Equation 40}
x^{2}-\left(x-ap^{1/2}+2\right)\left(x+ap^{1/2}-2\right)=\left(ap^{1/2}-2\right)^{2}.
\end{equation}

Note that increasing $x$ by $1$ in (\ref{Equation 7}) increases
the value of the expression by $2k+2x+2$.  Combining this with
(\ref{Equation 13}) and (\ref{Equation 40}), we find that one of
the quadratic residues in the aforementioned collection is
congruent (mod $p$) to an integer in the interval,
($\frac{p+1}{2}+k$, $\frac{p+1}{2}+k+\lfloor{p^{1/2}}\rfloor$) if

\[
\left(ap^{1/2}-2\right)^{2}>2k+2x+2-\lfloor p^{1/2}\rfloor,
\]
and since $x\leq{(1-a)p^{1/2}}$, we have

\begin{equation}\label{Equation 14}
\left(ap^{1/2}-2\right)^{2}>2k+2(1-a)p^{1/2}+2-\lfloor
p^{1/2}\rfloor,
\end{equation}
which proves the lemma.

\endproof

\begin{example}  {\em $p=13$ is the only prime number for which the greatest
number of consecutive quadratic non-residues modulo $p$ exceeds
$p^{1/2}$.}

\end{example}

\proof  Suppose $p>38659$, and suppose there exists a sequence of
more than $p^{1/2}$ consecutive quadratic non-residues (mod $p$).
As noted earlier, this implies that $p\equiv{13}$ (mod $24$).

Now suppose that

\begin{equation}\label{Equation 41}
\frac{p+1}{2}+k
\end{equation}
is a quadratic non-residue, where $k$ is a fixed integer. With
Lemma 1 in mind, we need only prove that this is not the least
quadratic non-residue in a sequence of more than $p^{1/2}$
consecutive quadratic non-residues when

\[
\frac{p^{1/2}}{2}+1<k<2^{1/2}p^{3/4}-p^{1/2}.
\]

Now consider three cases:

Case 1: $k<2p^{1/2}$.  Note that the difference between
(\ref{Equation 7}) at $x=\lfloor$$\frac{3p^{1/2}}{4}$$\rfloor$ and
$x=\lfloor$$\frac{p^{1/2}}{4}$$\rfloor$$-2$ is greater than $p$
because it equals

\[
\left((2k+1)\left\lfloor\frac{3p^{1/2}}{4}\right\rfloor+\left\lfloor\frac{3p^{1/2}}{4}\right\rfloor^{2}\right)-
\left((2k+1)\left(\left\lfloor\frac{p^{1/2}}{4}\right\rfloor-2\right)+\left(\left\lfloor\frac{p^{1/2}}{4}\right\rfloor-2\right)^{2}\right)
\]
\[
=2k\left(\left\lfloor\frac{3p^{1/2}}{4}\right\rfloor-\left\lfloor\frac{p^{1/2}}{4}\right\rfloor\right)+\left\lfloor\frac{3p^{1/2}}{4}\right\rfloor+4k
+\left\lfloor\frac{3p^{1/2}}{4}\right\rfloor^{2}-\left\lfloor\frac{p^{1/2}}{4}\right\rfloor^{2}+3\left\lfloor\frac{p^{1/2}}{4}\right\rfloor-2
\]
\[
>2k\left(\frac{p^{1/2}}{2}-1\right)+\frac{3p^{1/2}}{4}-1+4k+\frac{9p}{16}-\frac{3p^{1/2}}{2}+1-\frac{p}{16}+3\left(\frac{p^{1/2}}{4}-1\right)-2
\]
\[
=kp^{1/2}+2k+\frac{p}{2}-5>\frac{p}{2}+p^{1/2}+p^{1/2}+2+\frac{p}{2}-5>p,
\]
since $k>\frac{p^{1/2}}{2}+1$.  With Lemma 2 in mind, we find that
(\ref{Equation 41}) is not the least quadratic non-residue in a
sequence of more than $p^{1/2}$ consecutive quadratic non-residues
if

\begin{equation}\label{Equation 12}
\left(\frac{p^{1/2}}{4}-2\right)^{2}>2k+\frac{3p^{1/2}}{2}+2-\lfloor
p^{1/2}\rfloor.
\end{equation}
Since $k<2p^{1/2}$, (\ref{Equation 12}) holds whenever

\[
\left(\frac{p^{1/2}}{4}-2\right)^{2}>\frac{9p^{1/2}}{2}+3,
\]
which holds for $p>7711$.

Case 2: $2p^{1/2}<k<8p^{1/2}$.  Note that the difference between
(\ref{Equation 7}) at $x=\lfloor$$\frac{5p^{1/2}}{8}$$\rfloor$ and
$x=\lfloor$$\frac{3p^{1/2}}{8}$$\rfloor$$-2$ is greater than $p$
because it equals

\[
\left((2k+1)\left\lfloor\frac{5p^{1/2}}{8}\right\rfloor+\left\lfloor\frac{5p^{1/2}}{8}\right\rfloor^{2}\right)-
\left((2k+1)\left(\left\lfloor\frac{3p^{1/2}}{8}\right\rfloor-2\right)+\left(\left\lfloor\frac{3p^{1/2}}{8}\right\rfloor-2\right)^{2}\right)
\]
\[
=2k\left(\left\lfloor\frac{5p^{1/2}}{8}\right\rfloor-\left\lfloor\frac{3p^{1/2}}{8}\right\rfloor\right)+\left\lfloor\frac{5p^{1/2}}{8}\right\rfloor+4k
+\left\lfloor\frac{5p^{1/2}}{8}\right\rfloor^{2}-\left\lfloor\frac{3p^{1/2}}{8}\right\rfloor^{2}+3\left\lfloor\frac{3p^{1/2}}{8}\right\rfloor-2
\]
\[
>2k\left(\frac{p^{1/2}}{4}-1\right)+\frac{5p^{1/2}}{8}-1+4k+\frac{25p}{64}-\frac{5p^{1/2}}{4}+1-\frac{9p}{64}+3\left(\frac{3p^{1/2}}{8}-1\right)-2
\]
\[
=\frac{kp^{1/2}}{2}+2k+\frac{p}{4}+\frac{p^{1/2}}{2}-5>p+4p^{1/2}+\frac{p}{4}+\frac{p^{1/2}}{2}-5>p,
\]
since $k>2p^{1/2}$.  With Lemma 2 in mind, we find that
(\ref{Equation 41}) is not the least quadratic non-residue in a
sequence of more than $p^{1/2}$ consecutive quadratic non-residues
if

\begin{equation}\label{Equation 24}
\left(\frac{3p^{1/2}}{8}-2\right)^{2}>2k+\frac{5p^{1/2}}{4}+2-\lfloor
p^{1/2}\rfloor.
\end{equation}
Since $k<8p^{1/2}$, (\ref{Equation 24}) holds when

\[
\left(\frac{3p^{1/2}}{8}-2\right)^{2}>\frac{65p^{1/2}}{4}+3,
\]
which holds for $p>15917$.

Case 3: $8p^{1/2}<k$.  Note that the difference between
(\ref{Equation 7}) at $x=\lfloor$$\frac{17p^{1/2}}{32}$$\rfloor$
and $x=\lfloor$$\frac{15p^{1/2}}{32}$$\rfloor$$-2$ is greater than
$p$ because it equals

\[
\left((2k+1)\left\lfloor\frac{17p^{1/2}}{32}\right\rfloor+\left\lfloor\frac{17p^{1/2}}{32}\right\rfloor^{2}\right)-
\left((2k+1)\left(\left\lfloor\frac{15p^{1/2}}{32}\right\rfloor-2\right)+\left(\left\lfloor\frac{15p^{1/2}}{32}\right\rfloor-2\right)^{2}\right)
\]
\[
=2k\left(\left\lfloor\frac{17p^{1/2}}{32}\right\rfloor-\left\lfloor\frac{15p^{1/2}}{32}\right\rfloor\right)+\left\lfloor\frac{17p^{1/2}}{32}\right\rfloor+4k
+\left\lfloor\frac{17p^{1/2}}{32}\right\rfloor^{2}-\left\lfloor\frac{15p^{1/2}}{32}\right\rfloor^{2}+3\left\lfloor\frac{15p^{1/2}}{32}\right\rfloor-2
\]
\[
>2k\left(\frac{p^{1/2}}{16}-1\right)+\frac{17p^{1/2}}{32}-1+4k+\frac{289p}{1024}-\frac{17p^{1/2}}{16}+1-\frac{225p}{1024}+3\left(\frac{15p^{1/2}}{32}-1\right)-2
\]
\[
=\frac{kp^{1/2}}{8}+2k+\frac{p}{16}+\frac{7p^{1/2}}{8}-5>p+16p^{1/2}+\frac{p}{16}+\frac{7p^{1/2}}{8}-5>p,
\]
since $k>8p^{1/2}$.  With Lemma 2 in mind, we find that
(\ref{Equation 41}) is not the least quadratic non-residue in a
sequence of more than $p^{1/2}$ consecutive quadratic non-residues
if

\begin{equation}\label{Equation 30}
\left(\frac{15p^{1/2}}{32}-2\right)^{2}>2k+\frac{17p^{1/2}}{16}+2-\lfloor
p^{1/2}\rfloor.
\end{equation}
Since $k<2^{1/2}p^{3/4}-p^{1/2}$, (\ref{Equation 30}) holds
whenever

\[
\left(\frac{15p^{1/2}}{32}-2\right)^{2}>2^{3/2}p^{3/4}-\frac{31p^{1/2}}{16}+3,
\]
which holds for $p>27250$.

So when $p>38659$, no sequence of more than $p^{1/2}$ consecutive
quadratic non-residues exists.

Now all that remains is to consider the case $p\leq38659$.  This
case can be handled by a simple computation.  I have run a
computer program which compares the largest number of consecutive
quadratic non-residues modulo $p$ with $p^{1/2}$ for all primes
$p$, such that $p\equiv 13$ (mod $24$) and $p\leq38659$. From this
I was able to check that $13$ is the only prime number for which
the greatest number of consecutive quadratic non-residues (mod
$p$) exceeds $p^{1/2}$.

\endproof

\begin{remark}
The data obtained from this program can be viewed by going to the
website http://www.math.caltech.edu/people/hummel.html.  A sample
of some of the data obtained from the program is given below.  The
numbers in each set represent $p$, the greatest number of
consecutive quadratic non-residues (mod $p$), and
$\lfloor$$p^{1/2}$$\rfloor$ in that order.  For all but the
smallest numbers, $p^{1/2}$ far exceeds the greatest number of
consecutive quadratic non-residues.

$\{13,4,3\}$, $\{757,8,27\}$, $\{3181,9,56\}$, $\{5869,9,76\}$,
$\{7237,10,85\}$, $\{9397,10,96\}$, $\{12037,11,109\}$,
$\{14389,12,119\}$, $\{16477,12,128\}$, $\{18517,13,136\}$,
$\{20509,13,143\}$, $\{22381,12,149\}$, $\{24061,13,155\}$,
$\{26029,13,161\}$, $\{28429,13,168\}$, $\{30469,14,174\}$,
$\{32749,15,180\}$, $\{34693,14,186\}$, $\{36709,15,191\}$,
$\{38653,15,196\}$.
\end{remark}

Patrick Hummel is affiliated with the Class of 2006 at California
Institute of Technology.

\end{document}